\documentclass[12pt,a4paper,twoside,english]{article}
\usepackage{latexsym,graphicx}
\usepackage{color}
\usepackage{amsmath,amssymb}
\usepackage{babel}
\usepackage{graphicx}
\usepackage[latin1]{inputenc}
\newtheorem{teorema}{Theorem}
\newtheorem{corollario}[teorema]{Corollary}
\newtheorem{lemma}[teorema]{Lemma}
\newtheorem{proposizione}[teorema]{Proposition}
\newtheorem{definizione}[teorema]{Definition}

\newtheorem{esempio}[teorema]{Example}
\newtheorem{df}[teorema]{Definition}

\newtheorem{prop}[teorema]{Proposition}

\newcommand{\iE}{\mbox{$\displaystyle{\int_{E}}$}}

\newcommand{\erre}{\mbox{$\mathbb{R}$}}
\newcommand{\enne}{\mbox{$\mathbb{N}$}}

 \newcommand{\bvi}{\bigvee_{i=1}^{\infty} a_{i, \varphi(i)}}
 \newcommand{\bvib}{\bigvee_{i=1}^{\infty} b_{i, \varphi(i)}}

\newcommand{\bvik}{\bigvee_{i=1}^{\infty} a_{i, \varphi(i+k)}^k}

\title{\bf \Large Set-valued Kurzweil-Henstock Integral in Riesz spaces}
\author{\bf  A. Boccuto--A. M. Minotti --A. R. Sambucini
\footnote{
\noindent Supported by GNAMPA of CNR and University of Perugia \newline
corresponding author: Anna Rita Sambucini \newline
Authors' Address:
Dipartimento di Matematica e Informatica, via Vanvitelli,1
 I-06123 PERUGIA (Italy)
email: boccuto@dmi.unipg.it, minotti.am@gmail.com, matears1@unipg.it
}}

\date{}
\begin{document}
\maketitle \pagestyle{myheadings} \markboth{\centerline{\small \rm A. Boccuto -- A. M. Minotti -- A. R. Sambucini }}{\centerline{\small \rm
Set-valued Kurzweil-Henstock Integral in Riesz spaces}}
\small
\begin{abstract} A multivalued integral in Riesz spaces is given
using the Kurzweil-Henstock integral construction. Some of its properties and
a comparison with the Aumann approach  are also  investigated.
\end{abstract}
\mbox{~} \vskip.1cm
\noindent
\normalsize

\section{Introduction}

Riesz spaces and set-valued integrals were studied by several authors and in
different settings and have become an important tool in many branches of
applied mathematics, in particular in Economic Theory (see \cite{AB}),
where the theory of Riesz spaces is applied to study Riesz
commodities
and Riesz price economies,
utility functions and equilibria in an exchange economy.
Riesz spaces and ordered vector spaces provide a natural
framework for any modern theory of integration: indeed in this context
there are some relations with It\^{o} and Stratonovich-type integrals. Moreover
Riesz spaces have been recently used in the setting of
stochastic processes (see e.g. \cite{LW2}): indeed they can be viewed as
functions with values in the space $L^0(T,\Sigma,\mu)$ of all measurable functions
with identification up to $\mu$-null sets, where the convergence involved
is the almost everywhere convergence, which is not generated by any topology.
For a literature, see for example \cite{lavorodue,bck,BOCRIE,BRV,hlinena,pap,RN,R5,stajner}.

Multivalued analysis is a very powerful
 tool in the study of several problems
in many areas of Mathematics. For example, we recall here
subdifferentials of convex functionals, Calculus of Variations, degree theory,
fixed points, set-valued random processes, optimal control theory,
game theory, Pareto optimization, and so on.

The set-valued integration comes from Aumann using Bochner selections, and
a very large number of generalizations of this approach were given. For example
there are many recent studies concerning integrals of the type of  Sugeno, Choquet, Pettis,
Birkhoff, Kurzweil-Henstock,
pseudo-integrals (see \cite{CASCALES, CASCALES2, DPM, GSS}).
Following the ideas in \cite{JK, BS2004, BS} we extend in two directions
the results obtained. We introduce here a new kind of set-valued integral for
Riesz space-valued multifunctions
using the theory given in  \cite{BOCRIE}:
$\Phi(F,E)$ is in some sense the set of limits
with respect to the $(D)$-convergence
of all Riemann sums drawn out from
$\gamma$-fine partitions.
Its main properties like convexity, closedness, boundedness are investigated.
We also prove that in the case of simple multifunctions this integral agrees
with the usual integral obtained with Kurzweil-Henstock-type selections in the Aumann approach.

\section{Preliminaries}
\begin{definizione}\label{6Dsequence}
\rm Let $R$ be a Riesz space. A bounded double sequence $(a_{i,j})_{i,j}$ in $R$ is called a
\em $(D)$-sequence \rm or \textit{regulator} if $a_{i,j} \geq a_{i,j+1}$ for all $i,j
\in \enne$ and $\displaystyle{ \bigwedge_{j=1}^{\infty} \,
a_{i,j}=0}$ for all $i \in \enne $.\\
A Riesz space $R$ is said to be \textit{Dedekind complete} iff
every nonempty subset of $R$, bounded from above (below), has a lattice supremum (infimum) in
$R$ denoted by $\bigvee$ (resp. $\bigwedge$).
A Dedekind complete Riesz space $R$ is said to be \textit{super Dedekind complete} if every non-empty subset
$R_1 \subset R$,   $R_1 \neq \emptyset$,
bounded from above, contains a countable subset having the same supremum as $R_1$.

A Dedekind complete Riesz space $R$ is said to be \textit {weakly
$\sigma$-di\-stri\-bu\-ti\-ve} if for every
$(D)$-sequence $(a_{i,j})_{i,j}$ in $R$ one has:
\begin{eqnarray}{\label{6weaksigma}}
{\displaystyle \bigwedge_{\varphi \in
\enne^{\mathbb{N}}}\left(\bigvee_{i=1}^{\infty} \, a_{i,\varphi(i)}
\right)=0.}
\end{eqnarray}
\end{definizione}

We say that  $(R_1,R_2,R)$ is a  \em product triple \rm if there
exists a map $\cdot :R_1 \times R_2 \rightarrow R,$ which is compatible with the operations of sum
and order, suprema and infima (see also \cite[Assumption 5.2.1]{RN}).\\
 A Dedekind complete Riesz space $R$ is called an \em  algebra \rm if
 $(R,R,R)$ is a product triple. Our results are given in $(R,R,R)$ only for simplicity, but they hold
in general in a product triple $(R_1,R_2,R)$.\\
We assume that $(T,d)$ is a compact metric space, $\Sigma$ is the $\sigma$-algebra of all Borel
subsets of $T$ and $\mu: \Sigma \rightarrow R$ is a
positive additive mapping which is regular, namely:
for every $E \in \Sigma$ there exists a $(D)$-sequence
$(a^E_{i,j})_{i,j}$ such that for every $\varphi:\enne \rightarrow
\enne$ there exist  a compact set $K$  and  an open set $U$ such that $K \subset
E \subset U$ and
\begin{eqnarray}\label{reg}
\mu(U\setminus K) \leq \bigvee_{i=1}^{\infty} \, a^E_{i,\varphi(i)}.
\end{eqnarray}

We shall see in a moment that under usual conditions a regular measure is also $\sigma$-additive.

\begin{df}\label{6gage} \rm
A \em gage \rm on $T$ is a map $\gamma: T \rightarrow \mathbb{R}^+$.
\end{df}

\begin{df} \rm
A collection $\Pi$  of $T$ is a finite family
$\Pi = \{ (E_i,t_i): i=1, \ldots,k \} $ of couples such that
$t_i \in E_i, \, E_i \in \Sigma$ and $\mu(E_i\cap E_j) = 0$ for $i \neq j$.
 The points $t_i$, $i=1, \ldots, k$, are called \em
tags. \rm
 If moreover $\bigcup_{i=1}^k \, E_i= T$, $\Pi$ is called a \em partition. \rm
Given a gage $\gamma$, we say that $\Pi$ is \em $\gamma$-fine \rm $(\Pi \prec \gamma)$, if
$d(w,t_i) < \gamma (t_i)$
for every $i = 1 \ldots, k$\,
 and $ w \in E_i$.
\end{df}

A sequence $(r_n)_n$ in $R$ {\it $(D)$-converges} to an element $r\in R$ if there exists a
$(D)$-sequence $(a_{i,j})_{i,j}$ in $R$, such that
for every\, $\varphi \in \enne^{\mathbb{N}}$  there
exists  an integer $n_0$ such that  $|r_n-r| \leq
\bigvee_{i=1}^{\infty} a_{i, \varphi(i)},$ for all $n\geq n_0$.
In this case, we write $(D)\lim_n r_n=r$.

Observe that in weakly $\sigma$-distributive spaces the $(D)$-limits are unique.\\
We now recall the well-known Fremlin Lemma (see also
\cite{FREMLINLEMMA}):
\begin{lemma}{\label{6fremlin}} Let ${(a_{i,j}^{(k)})_{i,j}}$, $k \in \enne$, be any
countable family of $(D)$-sequences. Then for each fixed element
$u\in R$, $u \geq 0$, there exists a $(D)$-sequence
$(a_{i,j})_{i,j}$ such that, for every $\varphi\in \enne^{\mathbb{N}}$,
one has $$ u\wedge \sum _{k=1}^{s}
\left(\bigvee_{i=1}^{\infty}a_{i,\varphi (i+k)}^{(k)}\right)\leq
\bvi \quad \mbox{for all \,\,} \, s \in \enne.$$
\end{lemma}

From now on, we shall always suppose that $R$ is a
 weakly $\sigma$-distributive Dedekind complete  algebra.

  \bigskip

 As announced before we can prove that
\begin{teorema}
In  a weakly $\sigma$-distributive Dedekind complete  algebra every regular measure is $\sigma$-additive.
\end{teorema}
{\bf Proof}: Since $\mu$ is non-negative, all we must prove is that,
 for every disjoint sequence $(A_k)_k$ in $\Sigma$, one has
 $$\mu(\bigcup_{k\in \enne}A_k)\leq \sum_{k=1}^{+\infty}\mu(A_k).$$
 To this aim, let us denote by $A$ the union of all $A_k$'s, and observe that, for each integer $k$ there exists a $D$-sequence $(a_{i,j}^k)_{i,j}$,
 such that, for every $\varphi\in \enne^{\enne}$ it is possible to find a compact set $C_k$ and an open set $U_k$ satisfying $C_k\subset A_k\subset U_k$ and
 $\mu(U_k\setminus C_k)\leq \bvik$. Also, setting $u=2\mu(\Omega)$, by lemma \ref{6fremlin} it is possible to find a $D$-sequence $(b_{i,j})_{i,j}$
 such that
 $$ u\wedge \sum _{k=1}^{N}
\left(\bigvee_{i=1}^{\infty}a_{i,\varphi (i+k)}^{(k)}\right)\leq
\bvib \quad \mbox{for all \,\,} \, N \in \enne.$$
Now, by weak $\sigma$-distributivity and regularity, it is easy to deduce that $\mu(A)=\sup\{\mu(C):C\subset A, C\ {\rm compact} \}.$  Therefore, it will suffice to prove that
$$\mu(C)\leq \sum_{k=1}^{+\infty}\mu(A_k),$$
for every compact set $C\subset A$. So, fix any compact set $C\subset A$, and any mapping $\varphi\in \enne^{\enne}$. In correspondence with $\varphi$, there exist open sets $U_k$, $k=1,2...,$ such that $A_k\subset U_k$ and $\mu(U_k)\leq \mu(A_k)+\bvik$ for all $k$. Since $C$ is compact, and the union of all $U_k$'s contains $C$, there exists an integer $N$ such that $C\subset \cup_{k\leq N} U_k$. Since $\mu$ is non-negative we deduce:
$$\mu(C)\leq \sum_{k=1}^N\mu(U_k)\leq \sum_{k=1}^N \mu(A_k)+\sum_{k=1}^N \bvik=\mu(\bigcup_{k=1}^N A_k)+\sum_{k=1}^N \bvik, $$
and so
$$\mu(C)-\mu(\bigcup_{k=1}^NA_k)\leq u\wedge \sum_{k=1}^N \bvik\leq \bvib.$$
Then we have
$$\mu(C)\leq \mu(\bigcup_{k=1}^NA_k)+\bvib=\sum_{k=1}^N\mu(A_k)+\bvib\leq \sum_{k=1}^{+\infty}\mu(A_k)+\bvib.$$
But the last inequality holds for every $\varphi\in \enne^{\enne}$; so, since $R$ is weakly $\sigma$-distributive, we finally get
$$\mu(C)\leq \sum_{k=1}^{+\infty}\mu(A_k),$$
as desired.

\bigskip

We now give our definition of Kurzweil-Henstock integrability.
\begin{definizione}\rm \cite[Definition 3.1]{BOCRIE}\label{6equivalenzariesz}
\rm A function $f:T\rightarrow R$ is \em $(KH)$-integrable \rm
(or, in short, \em integrable \rm) if there exist $I \in R$ and a
$(D)$-sequence $(a_{i,j})_{i,j}$ such that for all $\varphi \in
\enne^{\mathbb{N}}$ there exists a gage $\gamma$ such that
for every
$\gamma$-fine
partition
$\Pi=\{(E_i, t_i), i=1, \ldots, q \}$   of $T$ we have:
\begin{eqnarray}\label{6trestelleriesz}
\left| \sum_{\Pi} \, f - I \right| \leq \bigvee_{i=1}^{\infty} \,
a_{i, \varphi(i)},
\end{eqnarray}
where  $\sum_{\Pi} \, f:= \sum_{i=1}^{q} \, f(t_i) \, \mu(E_i)$ is a Riemann sum.
The number $I$ is determined uniquely. It will be
denoted here simply by  $\int_T \, f \, d\mu.$
In an analogous way we can define the $(KH)$-integral on a set $E \in \Sigma$.
\end{definizione}

Observe that the $(KH)$-integral is additive on disjoint sets and is a
linear positive functional (see \cite[Propositions 3.2 and 3.5]{BOCRIE}). Moreover by
\cite[Lemma 1.10]{cinque} the integrability on $\Sigma$-measurable subsets is inherited.\\

In Riesz spaces in general it is not possible to compare Bochner and
Kurzweil-Henstock integrability,  we only know that
simple functions are integrable in both senses (\cite[Theorem 3.7]{BOCRIE}); moreover
it is possible to construct a Bochner integrable function which is not $(KH)$-integrable,
as the following example shows.

\begin{esempio}\label{7buoniebelli}
\rm Let $R_2 = \mathbb{R}$ and
$R  = R_1= c_{00}$ be the space of eventually null real-valued sequences. Let
$(u_n)_n$ be defined by
$u_n:=(0,\ldots, 0, 1, 0, \ldots)$,
where the value $1$ is assumed at the $n$-th coordinate.
The function $f:[0,1]\rightarrow R$, defined by
\begin{eqnarray}\label{7controesempio}
f(x)=\left \{
\begin{array}{ll}
u_n & {\rm if\,} x=1/n\\
\\0 & {\rm otherwise}
\end{array}\right.
\end{eqnarray}
vanishes
almost everywhere (with respect to the Lebesgue measure), so has null Bochner integral,
but  is not $(KH)$-integrable on $[0,1]$.

Indeed, fix arbitrarily $\delta:[0,1] \rightarrow
\erre^+$ and $n \in \enne$, $n \geq 2$. For every $i=1, \ldots,
n-1$, let $\displaystyle{\xi_i=\frac{1}{n+1-i}}$ and choose an
interval $]y_i,x_i[$ such that $\xi_i \in ]y_i, x_i[$, $x_i - y_i
< \delta(\xi_i)$, $ [y_i, x_i] \cap [y_j, x_j]= \emptyset$
for all $i \neq j$, $0 < y_1$ and $x_{n-1}<1$. We have:
$$0 <y_1 <x_1 < y_2 < x_2 < \ldots < y_{n-1} < x_{n-1} < 1.$$
Let $x_0=0$, $y_n=1$, and let us divide each of the intervals
$[x_{i-1}, y_i]$, $i=1, \ldots, n$, in subintervals, in such a way
to have a $\delta$-fine partition: this is possible, by virtue of
the Cousin Lemma (\cite[Theorem 2.3.1]{LEEV}).
These subintervals and the elements $([y_i,x_i],
\xi_i)$, $i=1, \ldots, n-1$, form a $\delta$-fine partition $\{
([t_{j-1}, t_j], \eta_j) :j=1, \ldots, p \}$. Since $f=0$ on each
of the intervals $[x_{i-1}, y_i]$, $i=1, \ldots, n$, we have:
$$\sum_{j=1}^{p} \,  (t_j - t_{j-1}) \, f(\eta_j) = \sum_{i=1}^{n-1}
\, (x_i - y_i) \, f(\xi_i).$$ Let $\lambda^{(n)}_i = x_{n+1-i} -
y_{n+1-i}$, $i=2, \ldots, n$: then we get
\begin{eqnarray}\label{7arrembaggio}
& & \sum_{j=1}^{p} \,  (t_j - t_{j-1}) \, f(\eta_j) = \nonumber
\sum_{i=1}^{n-1} \, \lambda^{(n)}_{n+1-i} \, f(\xi_i) =
\sum_{i=1}^{n-1} \, \lambda^{(n)}_{n+1-i} \, u_{n+1-i} \geq
\lambda^{(n)}_n \, u_n.
\end{eqnarray}
Since  $\lambda^{(n)}_n$ is
strictly positive for every $n$,
the sequence $\left(\sum_{\Pi(n)} \, f \right)_n$ is unbounded
in $R$.

If $f$ was $(KH)$-integrable on $[0,1]$, then there would exist a
gage $\delta_0:[0,1] \rightarrow \erre^+$ such that
\begin{eqnarray*}
\bigvee \left\{\sum_{\Pi} \, f : \Pi  {\rm \, is \, a \, }
\delta_0{\rm -fine \, partition \, of \, } [0,1] \right\} \in R :
\end{eqnarray*} this is a contradiction. Hence, $f$ is not
$(KH)$-integrable on $[0,1]$. $\quad \Box$
\end{esempio}


\section{Multivalued Kurzweil-Henstock integral}
We  introduce some definitions and properties.
\begin{definizione} \rm
Let ${\mathcal U}(C, r) := \{ z \in R: \exists \, x \in C :
|x-z| \leq r \}.$
A set $C \subset R$ is said to be \textit{closed} if $C = \mbox{cl}(C)$, where
\begin{eqnarray*}
\mbox{cl}(C) &:=& \bigcup_{(a_{i,j})_{i,j}}  \bigcap_{\varphi \in \enne^{\mathbb{N}}}
{\mathcal U}(C, \bigvee_{i=1}^{\infty} a_{i,\varphi(i)} )
\end{eqnarray*}
\end{definizione}
Let $cf(R)$ be the family of all non-empty, convex, bounded and closed subsets of a weakly $\sigma$-distributive Riesz space $R$. For example every interval   $[a,b]$ is closed. In fact, since
$\bigwedge_{\varphi} \bigvee_{i=1}^{\infty} a_{i,\varphi(i)} = 0$  for every ($D$)-sequence $(a_{i,j})_{i,j}$, then
 $$\mbox{cl}(C) := \bigcup_{(a_{i,j})_{i,j}}  \bigcap_{\varphi \in \enne^{\mathbb{N}}} [ a - \bigvee_{i=1}^{\infty} a_{i,\varphi(i)}, b + \bigvee_{i=1}^{\infty} a_{i,\varphi(i)}] = [a,b].$$
Analogously as in \cite{CV} the addition
$\stackrel{\cdot}{+}: cf(R) \times cf(R) \rightarrow cf(R)$ is defined as follows:
it  is the closure of  the direct sum of the two sets. If $A_i,\,\, i = 1, \ldots, n$ are in $cf(R)$ we denote by $\sum_{i=1}^n A_i$ the set
\begin{eqnarray}\label{sommapunto}
\sum_{i=1}^n A_i:= \mbox{cl}(A_1  +  \cdots  + A_n).
\end{eqnarray}
Observe that if $A,B \in cf(R)$ then $A \stackrel{\cdot}{+} B \in cf(R)$. We have only to check the convexity. So, let $z_1, z_2 \in A \stackrel{\cdot}{+} B$, we shall prove that  for every  real number $\alpha \in ]0,1[$,
$$\alpha z_1 + (1-\alpha) z_2 \in A \stackrel{\cdot}{+} B.$$
 Let $(a^{(1)}_{i,j})_{i,j}$, $(a^{(2)}_{i,j})_{i,j}$
be two regulators associated with $z_1$, $z_2$ respectively, and
set $b_{i,j}=2(a^{(1)}_{i,j}+a^{(2)}_{i,j})$ for all $i,j \in
\enne$. It is easy to check that the double sequence
$(b_{i,j})_{i,j}$ is a regulator. So for every $\varphi \in \enne^{\mathbb{N}}$ and $l=1,2$ there exist $x_l \in A, y_l  \in B$ such that
\[ | z_l - (x_l + y_l) | \leq \bigvee_{i=1}^{\infty} a^{(l)}_{i, \varphi(i)}.\]
Then
\begin{eqnarray*}
| \alpha z_1 + (1 - \alpha) z_2)| \leq  |\alpha( x_1 + y_1) + (1 - \alpha) (x_2 + y_2)| + \bigvee_{i=1}^{\infty} b_{i, \varphi(i)}
\end{eqnarray*}
and this proves the convexity, since the direct sum is convex too.

\begin{definizione}\label{stronglybounded} \rm
A multifunction $F$ is \it bounded, \rm if
there exists a positive element $L \in R$ with
$F(t) \subset
[-L,L]$ for all $t \in T$.
\end{definizione}

We now define a multivalued integral as follows:
\begin{definizione}\label{starintegral}
\rm Let $F:T \to 2^R$ be a multifunction, and $E \in \Sigma$.
We call \it $(*)$-integral \rm of $F$ on $E$ the set
\begin{eqnarray*}
\Phi(F,E) &=& \left\{ \right. z \in R: \text{ there exists } \,\, (a_{i,j})_{i,j}:
\text{ for all } \varphi \in \enne^{\mathbb{N}}\,\, \text { there is a gage } \, \gamma \\ &&
\text{such that for every\,\,}
\gamma\text{-fine partition }  P_{\gamma}:=
\{ (E_i,t_i): i=1, \ldots,k \} \in \Pi_{\gamma} \\
&&
\text{of } E \text{ there exists }
\, \displaystyle{ \, c \in
\sum_{i=1}^{k} \, F(t_i) \, \mu(E_i)} \text{ with\,\,}
\displaystyle{|z-c| \leq \bigvee_{i=1}^{\infty}
a_{i, \varphi(i)}}\,
\left. \right\}.
\end{eqnarray*}
\end{definizione}

It is easy to check that, if $F$ is single-valued and integrable, then its
$(*)$-integral coincides with the usual Kurzweil-Henstock
integral given in \cite{BOCRIE}.
Observe that, in the case of Example \ref{7buoniebelli}, we get $\Phi(\{f\},E) = \emptyset$ since the
Riemann sums are not bounded in $R$.
\\

Then we have the following characterization:
\begin{prop}\label{estrella}
The  set $\Phi(F,E)$  can be obtained as follows:
$$ \Phi(F,E)=\bigcup_{(a_{i,j})_{i,j}} \, \bigcap_{\varphi \in \enne^{\mathbb{N}}}
\bigcup_{\gamma} \bigcap_{\{ (E_i,t_i)\} \in \Pi_{\gamma}} \,
{\mathcal U} \left( \sum_{i=1}^{k}\, F(t_i) \mu(E_i),
\bigvee_{i=1}^{\infty} a_{i,\varphi(i)} \right).$$
Moreover, if $u$ is an upper bound for $F$, then $u\, \mu(T)$ is an upper bound for $\Phi(F,E)$ for every $E \in \Sigma$.
\end{prop}
{\bf Proof}:
This  is an easy consequence of definitions of
$(*)$-integral and of $\bigvee_{i=1}^{\infty} a_{i,\varphi(i)}$-neigh\-bor\-hood.
For the second part,
thanks to boundedness and weak $\sigma$-distributivity, we have:
\begin{eqnarray*}
\Phi(F,E)
&\subset&
\bigcup_{(a_{i,j})_{i,j}} \, \bigcap_{\varphi \in \enne^{\mathbb{N}}}
\bigcup_{\gamma} \bigcap_{\{ (E_i,t_i)\} \in \Pi_{\gamma}} \,
{\mathcal U} \left( [-u \mu(T), u \mu(T)],
\bigvee_{i=1}^{\infty} a_{i,\varphi(i)} \right) \subset
\\ &\subset&
\bigcup_{(a_{i,j})_{i,j}} \,\,
\bigcap_{\varphi \in \enne^{\mathbb{N}}} [-u \mu(T) -
\bigvee_{i=1}^{\infty} a_{i,\varphi(i)}, u \mu(T) + \bigvee_{i=1}^{\infty} a_{i,\varphi(i)}]
= [-u \mu(T), u \mu(T)].
\end{eqnarray*}
$\Box$\\

We now prove  convexity and  closedness of the set $\Phi(F,E)$.
\begin{prop}
$\Phi(F,E)$ is convex  provided that it is non-empty and the multifunction $F$ is convex-valued.
\end{prop}
{\bf Proof}:
 Let $z_1$, $z_2
\in \Phi(F,E)$; $(a^{(1)}_{i,j})_{i,j}$, $(a^{(2)}_{i,j})_{i,j}$
be two regulators associated with $z_1$, $z_2$ respectively, and
set  $b_{i,j}=2(a^{(1)}_{i,j}+a^{(2)}_{i,j})$ for all $i,j \in
\enne$.  Since $z_1$ and $z_2$ belong to
$\Phi(F,E)$, then to every mapping $\varphi \in \enne^{\mathbb{N}}$
there correspond two gages $\gamma_1$, $\gamma_2$ such
 that, whenever $\Pi:=\{(G_l,\eta_l):l=1, \ldots, q\}$ is
both a $\gamma_1$- and a $\gamma_2$-fine partition, two points
$\displaystyle{c_1, c_2 \in \sum_{l=1}^q F(\eta_l) \, \mu(G_l)}$
can be found, with
$$|z_s-c_s| \leq \bigvee_{i=1}^{\infty} a^{(s)}_{i,\varphi(i)}, \quad s=1,2.$$
Let now $\gamma(x):=\gamma_1(x) \wedge \gamma_2(x)$, $x \in T$. Of
course, every $\gamma$-fine partition is both $\gamma_1$- and
$\gamma_2$-fine. Moreover, note that the set
$\displaystyle{\sum_{l=1}^q F(\eta_l) \, \mu(G_l)}$ is convex,
since $F$ is convex-valued.
Hence, for every $\alpha \in [0,1]$,
$\displaystyle{\alpha c_1 + (1-\alpha) \, c_2 \in \sum_{l=1}^q
F(\eta_l) \, \mu(G_l)}$, and

\begin{eqnarray*}
|\alpha \, c_1 +(1-\alpha) c_2 - \alpha z_1 -(1-\alpha) z_2|
&\leq&
 \alpha |c_1 - z_1 | + (1-\alpha) |c_2 - z_2 | \leq\\
&\leq &
\alpha \bigvee_{i=1}^{\infty} b_{i,\varphi(i)} + (1-\alpha)
\bigvee_{i=1}^{\infty} b_{i,\varphi(i)}=\bigvee_{i=1}^{\infty}
b_{i,\varphi(i)}.
\end{eqnarray*}
 Thus, $\alpha \, z_1 + (1-\alpha) \, z_2 \in
\Phi(F,E)$, which proves convexity of $\Phi(F,E)$. $\quad \Box$ \vspace{3mm}
\\
\begin{prop}\label{chiusura}

If $F$ is bounded and $R$ is also super Dedekind complete, then the set $\Phi(F,E)$ is closed.
\end{prop}
{\bf Proof}:
We have to prove that
\[ \Phi(F,E) = \bigcup_{(a_{i,j})_{i,j}} \, \bigcap_{\varphi \in \enne^{\mathbb{N}}}
{\mathcal U}\Bigl(\Phi(F,E), \bigvee_{i=1}^{\infty} a_{i,\varphi(i)}\Bigr).\]
One
 inclusion is obvious, for the converse let
$z \in \bigcup_{(a_{i,j})_{i,j}} \, \bigcap_{\varphi \in \enne^{\mathbb{N}}}
{\mathcal U}\Bigl(\Phi(F,E), \bigvee_{i=1}^{\infty} a_{i,\varphi(i)}\Bigr)$.
Then there exists a regulator $(a_{i,j})_{i,j}$ such that for every
$\varphi \in \enne^{\mathbb{N}}$ there exists $c_{\varphi} \in \Phi(F,E)$ with
\[ | z - c_{\varphi} | \leq \bigvee_{i=1}^{\infty} a_{i,\varphi(i)}.\]
Since $R$ is super Dedekind complete and weakly $\sigma$-distributive,
a sequence $(\varphi_n)_n \in \enne^{\mathbb{N}}$  can be found such that
\[ 0 = \bigwedge_{\varphi}\left(\bigvee_{i=1}^{\infty} a_{i,\varphi(i)} \right)=
\bigwedge_{n} \left(\bigvee_{i=1}^{\infty} a_{i,\varphi_n(i)} \right).\]
Without loss of generality, arguing as in \cite{Sobham}, we can assume that
$(\varphi_n)_n$ is increasing and $\varphi_n(i) \leq \varphi_n (i+1)$ for every $i,n \in \enne$.
Let now $z_n = c_{\varphi_n}$, so that
 $z =
(o)\lim_n \, z_n =
(D)\lim_n \, z_n$, where $z=(o)\lim_n \, z_n$ means that $\wedge_n \sup_{m \geq n} |z_m - z|= 0$.
\\
We will show that  $z \in \Phi(F,E)$. Observe that the sequence $(z_n)_n$ is
bounded. Since $z_n \in \Phi(F,E)$, a  regulator
$(\alpha^{(n)}_{i,j})_{i,j}$ can be found, with the property that to
every $\varphi \in \enne^{\mathbb{N}}$ there corresponds a gage
$\gamma_n$ such that for any $\gamma_n$-fine partition $\{
(E_i,t_i):i=1, \ldots,k \}$ there is $w_n \in \sum_{i=1}^k F(t_i)
\, \mu(E_i)$ with $$|z_n - w_n | \leq \bigvee_{i=1}^{\infty}
\alpha^{(n)}_{i, \varphi(i+n)}$$ for all $n \in \enne$. Since $F$ is
bounded, then $(w_n)_n$ is bounded, and hence the
sequence $(z_n - w_n)_n$ is bounded too. Thus, by virtue of Lemma \ref{6fremlin},
there is a regulator $(\beta_{i,j})_{i,j}$ with
$$|z_n - w_n | \leq \bigvee_{i=1}^{\infty} \, \beta_{i, \varphi(i)}$$
for any $n \in \enne$. Moreover, since the sequence $(z_n)_n$
$(D)$-converges to $z$,  a regulator $(b_{i,j})_{i,j}$ can be
found, with the property that to every $\varphi \in \enne^{\mathbb{N}}$
there corresponds $\overline{n} \in \enne$ such that
$$|z_n - z| \leq \bigvee_{i=1}^{\infty} \, b_{i,\varphi(i)} \quad {\rm whenever \, \,}
n \geq \overline{n}.$$
Let now
$c_{i,j}=2(\beta_{i,j} + b_{i,j})$, $i,j \in \enne$. It is easy to see
that $(c_{i,j})_{i,j}$ is a regulator.\\
The element $z_{\overline{n}}$ is such that
for every $\varphi \in \enne^{\mathbb{N}}$ there is a gage
$\gamma_{\overline{n}}$ with the property that for all
$\gamma_{\overline{n}}$-fine partition $\{ (E_i,t_i):i=1, \ldots,k
\}$ there is $w_{\overline{n}} \in \sum_{i=1}^k F(t_i) \,
\mu(E_i)$ with
$$|z_{\overline{n}} - w_{\overline{n}} | \leq \bigvee_{i=1}^{\infty} \beta_{i, \varphi(i)}$$
and
$$|z - w_{\overline{n}}| \leq|z -z_{\overline{n}}|+ |z_{\overline{n}} -
w_{\overline{n}}| \leq \bigvee_{i=1}^{\infty} c_{i, \varphi(i)}.$$
This implies that $z \in \Phi(F,E)$, and hence the set $\Phi(F,E)$
is closed. $\quad \Box$\\

\begin{proposizione}
If  $0 \in F(t)$ for every $t \in T$, then
\[ \Phi(F,A) \subset \Phi(F,B) \,\, \mbox{for every} \,\, A, B \in \Sigma, \, A \subset B.\]
\end{proposizione}
{\bf Proof}:
Every element of the Riemann sum $\sum_{\Pi} F|_A$ can be extended
to the Riemann sum $\sum_{\Pi} F|_B$, using \cite[Lemma 1.2]{cinque} (Cousin Lemma)
and the fact that $0 \in F(t_i)$, for every $t_i \in B \setminus A$. Now
the proof follows using the definition of $(*)$-integral. $\Box$\\

The next goal, that we will  obtain in the next section, is to prove that simple measurable multifunctions are
integrable and the integrals can be obtained as in the single-valued case.
We prove this fact in several steps.
We begin with the following
\begin{proposizione}\label{costante}
Let $C$ be a closed subset of $R$. If  $F(t) \equiv C$ then, for every $E \in \Sigma$, we get
\begin{eqnarray*}
\Phi(F,E) = \Phi(F 1_E, T) = C \mu(E).
\end{eqnarray*}
\end{proposizione}
{\bf Proof:}
We will  show that
\[ C\mu(E) \subset \Phi(C 1_E,T) \subset \Phi(C,E) \subset C \mu(E).\]
First of all observe that, since $\mu$ is regular,
in correspondence with $E$
there is a regulator $(a^E_{i,j})_{i,j}$ such that for every $\varphi: \enne \rightarrow \enne$,
 there exist a compact set $K$ and an open set $U$
such that $K \subset E \subset U$ and satisfying (\ref{reg}).
Arguing analogously
as in \cite{BOCRIE}, since $K$ is compact and $U$ is open, there exists a gage $\gamma^E$ such
that
\[B(t,\gamma^E(t)) :=  \{ w \in T: d(t,w) < \gamma^E (t) \} \subset
\left\{ \begin{array}{ll}
  U & \forall\,\, t \in K,\\
  U \setminus K & \forall\,\, t \in U \setminus K,\\
 T \setminus K  & \forall\,\, t \not \in U.
\end{array} \right. \]
So, if
$\Pi= \{ (D_i, u_i), i=1,\ldots, q \}$ is a $\gamma^E$-fine partition of $T$, then
$$\mu( E \Delta (\cup \{D_i :  u_i \in E \})) \leq \bigvee_{i=1}^{\infty} a^E_{i, \varphi(i)}.$$
Let $z \in C \mu(E)$ and consider  $(a^E_{i,j})_{i,j}$ and $\gamma^E$ as above.
Let $\Pi= \{ (D_i, u_i), i=1,\ldots, q \}$ be  a $\gamma^E$-fine partition of $T$, and  $\widetilde{D} =
\cup_i \{D_i : u_i \in E\}$. Then
\[ \sum_{i=1}^q C 1_E (u_i) \mu(D_i) = \sum_{u_i \in E} C 1_E (u_i) \mu(D_i) +
	\sum_{u_i \not\in E} C 1_E (u_i) \mu(D_i) =
C \mu(\widetilde{D}).\]
Since $\mu(\widetilde{D}) \in [ \mu(E) - \bigvee_{i=1}^{\infty}
a^E_{i,\varphi(i)},  \mu(E) + \bigvee_{i=1}^{\infty} a^E_{i,\varphi(i)}]$, then there exists $c \in  \sum_{i=1}^q C 1_E (u_i) \mu(D_i)$ such that $ |z - c| \leq \bigvee_{i=1}^{\infty} a^E_{i,\varphi(i)}$
and this proves the first inclusion. \\
If  $z \in \Phi(F 1_E, T)$, then there is a regulator
$(b_{i,j})_{i,j}$ such that for every  $\varphi \in \enne^{\mathbb{N}}$,
there exists a gage $\gamma$ such that for every
$\gamma$-fine partition $\Pi = \{ (T_i,t_i):i=1, \ldots, q\} $
there is $c \in \sum_{\Pi} F 1_E$ such that
$|z -c| \leq \bigvee_{i=1}^{\infty} b_{i,\varphi(i)}$ .  In particular this holds for  $\gamma$-fine partitions $\Pi^{\prime}= \{(D_i,t_i): i =1, \ldots,s\}$
such that $\{(D_i,t_i):t_i \in E\} $ are partitions of $E$.  This proves that $z \in \Phi(F,E)$.\\
For the last inclusion,
 if we take $z \in \Phi(F,E)$, then there exists a regulator
$(a_{i,j})_{i,j}$ such that for every  $\varphi \in \enne^{\mathbb{N}}$ there is a gage
$\gamma$ such that for every $\gamma$-fine partition $\{ (E_i,t_i):i=1, \ldots, k\}
\in \Pi_{\gamma}$ of $E$ there
is $c \in\sum_{i=1}^{k} \, F(t_i) \, \mu(E_i)$  with $|z-c| \leq \bigvee_{i=1}^{\infty}
a_{i, \varphi(i)}$. \\
Since the partition  $\{ (E_i,t_i): i=1, \ldots, k\}$
is given for the set $E$, then the tags $t_i$ are in $E$ and so
we get $F(t_i) = C$ for every $i =1, 2, \ldots, k $.
Then $c \in C\mu(E)$ and $z\in {\mathcal U}(C\mu(E), \bigvee_{i=1}^{\infty}
a_{i, \varphi(i)})$.
Since this holds for every $\varphi$ and $C$ is closed,
we get $z\in C\mu(E)$ and finally $\Phi(F,E) \subset C\mu(E)$ by arbitrariness of $z$, and this proves the last inclusion. \\
 $\quad \Box$ \vspace{3mm}

\subsection{Comparison with the Aumann integral in Riesz spaces}
We introduce now the Aumann integral via
Kurzweil-Henstock integrable selections,
in order to compare it with the previous
integral for multifunctions. For a  multifunction
\mbox{$F: T \rightarrow cf(R)$}
let
\( S^1_F =\{  f: f(t) \in F(t)
\hskip.1cm \mu-\mbox{a.e. and \,}  f \,\,  \mbox{is $(KH)$-integrable } \}
\)
be the set of all $(KH)$-integrable selections of $F$.

\begin{definizione} \rm
If $F$ is such that $S_F^1$ is non-empty,
then for every $E \in \Sigma$ we define the
{\em Aumann integral~} (shortly $(A)$-{\em integral~}) as
\[ (A)\iE F d \mu = \left\{  \iE f d\mu, f \in S^1_F \right\}.\]
\end{definizione}

We recall
 that, in the context of Banach spaces, the Aumann integral is defined via Bochner integrable selections.
In our setting, this is not a good idea, as we showed in Example  \ref{7buoniebelli}.
As in the single-valued case we obtain that
\begin{teorema}\label{confronto}
If $F = \sum_{k=1}^n C_k 1_{E_k}$ is closed valued, $E_k \in \Sigma$ for
all $k=1$, $\ldots$, $n$ and the $E_k$'s are pairwise disjoint, then for every $A \in \Sigma$ we get
$$\mbox{cl}\left(\Phi(\sum_{k=1}^n C_k 1_{E_k}, A) \right)=
 \sum_{k=1}^n \Phi(C_k 1_{E_k}, A) =
\sum_{k=1}^n C_k \mu(E_k \cap A) =
\mbox{cl}\left( (A)\int_A F d \mu \right).$$
\end{teorema}
{\bf Proof}:
First of all observe that in this case all the sets involved are non-empty.
The equality
$ \sum_{k=1}^n \Phi(C_k 1_{E_k}, A) =
\sum_{k=1}^n C_k \mu(E_k \cap A)$
follows immediately from
Proposition \ref{costante}
with $T=A$, $E=E_k$, $C=C_k$, $k=1$, $\ldots$, $n$.
Moreover observe that, by (\ref{sommapunto}),  if $F = \sum_{k=1}^n C_k 1_{E_k}$, then for every gage $\gamma$ and  for all $\gamma$-partitions
$\Pi = \{ (B_r, t_r): r= 1, \ldots w \}$ of $A$ we have:
\begin{eqnarray}\label{commuta}
\sum_{\Pi} F &=&
\sum_{r=1}^w F(t_r) \mu (B_r) =
\sum_{r=1}^w \left( \sum_{k=1}^n C_k 1_{E_k} \right) (t_{r}) \mu(B_r) =\\ \nonumber &=&
 \sum_{k=1}^n C_k \left( \sum_{r=1}^w 1_{E_k} (t_{r}) \mu(B_r) \right) = \sum_{k=1}^n  \sum_{\Pi} F|_{E_k}.
\end{eqnarray}
We  now prove that
\begin{eqnarray}\label{mialc}
 \Phi(\sum_{k=1}^n C_k 1_{E_k}, A) \subset
 \sum_{k=1}^n \Phi(C_k 1_{E_k}, A) .
 \end{eqnarray}
Let $z \in \Phi(\sum_{k=1}^n C_k 1_{E_k}, A)$, then  there is a
regulator $(b_{i,j})_{i,j}$ with the property that for any $\varphi
\in \enne^{\mathbb{N}}$ there is a gage $\gamma$ such
that  for all $\gamma$-fine partitions
$\Pi_{*}=\{(D_m,\xi_m):m=1, \ldots, q \}$ of $A$ there is
$y \in  \sum_{\Pi_*} \left( \sum_{k=1}^n C_k 1_{E_k} \right) $
such that
$\displaystyle{|y-z|\leq \bigvee_{i=1}^{\infty} b_{i,\varphi(i)}}$.
We consider now only the $\gamma$-fine partitions
$\Pi^{*}=\{(D^{\prime}_m ,t_{m}): m=1, \ldots, s\}$ of $A$  such that for every
$k=1, \ldots, n$ the family
$\Pi^{*}|_{E_k}$  is also a partition of $E_k$.
By (\ref{commuta}), we have that
$$y \in \sum_{m=1}^s \left( \sum_{k=1}^n C_k 1_{E_k} \right) (t_{m}) \mu(D^{\prime}_m)
= \sum_{k=1}^n C_k \left( \sum_{m=1}^s 1_{E_k} (t_{m}) \mu(D_m) \right).$$
Thus for all $k=1, \ldots, n$ there are a regulator $(c_{i,j})_{i,j})$ and
$y_k \in   \Phi(C_k 1_{E_k}, A)$, for $k=1, \ldots, n$, such that for every $\varphi \in \enne^{\mathbb{N}}$ we have
 $\displaystyle{|y- \sum_{k=1}^n \, y_k}| \leq \bigvee_{i=1}^{\infty} c_{i, \varphi(i)}$.
This implies  that
\[z  \in   \bigcap_{\varphi \in \enne^{\mathbb{N}}}
{\mathcal U} \left(\sum_{k=1}^n \Phi(C_k 1_{E_k},A),
\bigvee_{i=1}^{\infty} 2( b_{i, \varphi(i)} + c_{i, \varphi(i)})
\right)
\subset \bigcup_{a_{i,j}} \bigcap_{\varphi \in \enne^{\mathbb{N}}}
{\mathcal U} \left(\sum_{k=1}^n \Phi(C_k 1_{E_k},A),
\bigvee_{i=1}^{\infty} a_{i, \varphi(i)} \right).\]
Since  $\sum_{k=1}^n \Phi(C_k 1_{E_k},A)$ is closed, then
 the inclusion in (\ref{mialc}) holds.
\\
 Obviously
$\sum_{k=1}^n \Phi(C_k 1_{E_k}, A)  = \sum_{k=1}^n C_k \mu(A \cap E_k) \subset (A)\int_A F d \mu \subseteq
\mbox{cl}\left( (A)\int_A F d \mu \right).$\\
So we have only to prove that
$(A)\int_A F d \mu \subset \Phi(F,A)$.
But, if $z \in (A)\int_A F d \mu$, then $z = \int_A g d\mu$ with
$g \in S^1_F$. By Definition \ref{6equivalenzariesz}, there exists
a regulator $(a_{i,j})_{i,j}$ such that for all $ \varphi \in
\enne^{\mathbb{N}}$ there exists a gage $\gamma$ such that
for every  $\gamma$-fine
partition
$\Pi$ we have:
\begin{eqnarray}
\left| \sum_{\Pi} \, g - z \right| \leq \bigvee_{i=1}^{\infty} \,
a_{i, \varphi(i)}.
\end{eqnarray} Since $\sum_{\Pi} \, g  \in \sum_{\Pi} \, F$, then the assertion follows.
$\Box$\\

\begin{corollario}
Suppose that $R$ is super Dedekind complete. If $F$ is a simple multifunction with closed and bounded values then
$$\Phi(\sum_{k=1}^n C_k 1_{E_k}, A)  =
\mbox{cl}\left( (A)\int_A F d \mu \right).$$
\end{corollario}
{\bf Proof}: it is an immediate consequence of Theorem \ref{confronto}  since $\Phi(\sum_{k=1}^n C_k 1_{E_k}, A)$ is closed thanks to Proposition \ref{chiusura}. $\Box$\\

It is still an open problem to compare the (*)- and Aumann integrals for not necessarily simple multivalued functions.

\end{document}